\documentclass{article}
\pdfoutput=1

\usepackage[margin=1in]{geometry}
\usepackage[utf8]{inputenc}
\usepackage[algoruled,linesnumbered]{algorithm2e}
\usepackage{amsfonts}       
\usepackage{amsmath}
\usepackage{amssymb}
\usepackage{booktabs}       
\usepackage[font=small]{caption}
\usepackage{doi}
\usepackage{enumitem}
\usepackage{float}
\usepackage{graphicx}
\usepackage{hyperref}       
\usepackage{microtype}      
\usepackage{multirow}
\usepackage{nicefrac}       
\usepackage{soul}
\usepackage{subcaption}
\usepackage{url}            
\usepackage{xcolor}

\providecommand{\keywords}[1]
{
  \hspace{2.5em}\small
  \textbf{Keywords:} #1
}
\newcommand{\abs}[1]{\left\lvert#1\right\rvert}
\newcommand{\cardinality}[1]{\left\lvert#1\right\rvert}
\newcommand{\norm}[1]{\left\lVert#1\right\rVert}

\newcommand{\redtext}[1]{{\color{black}{#1}}}
\newcommand{\greentext}[1]{{\color{black}{#1}}}
\newcommand{\bluetext}[1]{{\color{black}{#1}}}
\newcommand{\yellowtext}[1]{{\color{black}{#1}}}
\newcommand{\orangetext}[1]{{\color{black}{#1}}}
\newcommand{\tp}[2][-9]{{#2}^{\mkern#1mu\top}} 

\title{\textbf{Computing the k-th Eigenvalue of Symmetric $H^2$-Matrices}}

\author{M. Ridwan Apriansyah$^1$ and Rio Yokota$^2$}
\date{%
    $^1$School of Computing, Tokyo Institute of Technology\\%
    \texttt{ridwan@rio.gsic.titech.ac.jp}\\%
    $^2$Global Scientific Information and Computing Center, Tokyo Institute of Technology\\%
    \texttt{rioyokota@gsic.titech.ac.jp}\\
}

\hypersetup{
pdftitle={Computing the k-th Eigenvalue of Symmetric $H^2$-Matrices},
pdfsubject={},
pdfauthor={M. Ridwan Apriansyah and Rio Yokota},
pdfkeywords={eigenvalue solver, bisection, $H^2$-matrix, generalized LDL factorization},
}

\begin{document}
\maketitle

\begin{abstract}
The numerical solution of eigenvalue problems is essential in various application areas of scientific and engineering domains. In many problem classes, the practical interest is only a small subset of eigenvalues so it is unnecessary to compute all of the eigenvalues. Notable examples are the electronic structure problems where the $k$-th smallest eigenvalue is closely related to the electronic properties of materials. In this paper, we consider the $k$-th eigenvalue problems of symmetric dense matrices with low-rank off-diagonal blocks. We present a linear time generalized LDL decomposition of $\mathcal{H}^2$ matrices and combine it with the bisection eigenvalue algorithm to compute the $k$-th eigenvalue with controllable accuracy. In addition, if more than one eigenvalue is required, some of the previous computations can be reused to compute the other eigenvalues in parallel. Numerical experiments show that our method is more efficient than the state-of-the-art dense eigenvalue solver in LAPACK/ScaLAPACK and ELPA. Furthermore, tests on electronic state calculations of carbon nanomaterials demonstrate that our method outperforms the existing HSS-based bisection eigenvalue algorithm on 3D problems.
\end{abstract}
\section{Introduction}
\label{section:introduction}
We consider the standard eigenvalue problem of the form
\begin{equation}
    Av = \lambda v,
\end{equation}
where $A$ is $n \times n$ real symmetric matrix. This arises in a wide range of applications in science and engineering, ranging from electronic structure calculations to vibration analysis and structural dynamics~\cite{Andrew2013}.
Commonly used solutions rely on tridiagonal reduction and iterative solvers based on QR iteration, bisection, and divide-and-conquer algorithms~\cite{Trefethen1997}. For general symmetric matrix $A$, this typically requires $O(n^3)$ operations.

In recent decades, many new dense eigenvalue solvers have emerged, where the low-rank structure of the matrix is exploited to accelerate the computations. In this work, we consider such class of matrices with many rank-deficient off-diagonal blocks. This kind of matrix is often encountered in computational science, notably from the discretization of elliptic partial differential operators that govern a wide range of application areas~\cite{Hackbusch2015}. 
There are several kinds of structured low-rank formats, such as BLR~\cite{Amestoy2015_BLR}, BLR\textsuperscript{2}~\cite{Ashcraft2021_BLR2}, HODLR~\cite{Ambikasaran2013_HODLR}, quasiseparable/semiseparable~\cite{Chandrasekaran2005_SSS,Vandebril2008_QSS}, HSS~\cite{Chandrasekaran2007_HSS}, $\mathcal{H}$~\cite{Hackbusch1999_H}, and $\mathcal{H}^2$-matrices~\cite{Hackbusch2002_H2}.
Many studies have been conducted to obtain fast algorithms by adapting well-known dense eigenvalue solvers to these formats. QR iterations have been studied with quasiseparable matrices to compute all eigenvalues of generalized companion matrices in $\mathcal{O}(n^2)$ operations~\cite{Chandrasekaran2008_QR_QS,VanBarel2010_QR_QS}.
A similar iterative method based on LR Cholesky transformations has also been studied with HODLR matrices~\cite{Benner2013_LRCholesky}.
Further, divide-and-conquer eigenvalue solvers have been adapted to HODLR and HSS matrices, reducing the cost to compute all eigenvalues down to $\mathcal{O}(n\log_2^2(n))$ operations~\cite{Ou2022_HSS_DnC,Susnjara2021_HODLR_DnC}.

However, in some problem classes, the practical interest is only a small subset of eigenvalues so it is unnecessary to compute all of the eigenvalues~\cite{Morgan1991_InteriorEigenvalue}. Notable examples are the electronic structure calculations of materials~\cite{Hoshi2012_OrderNElectronicStructure,Lee2018_KthEigenvalue}, weather forecasting~\cite{Tanaka1985_WeatherForecasting}, and study of tidal motion~\cite{Cline1976_TidalMotion}. In such cases, the existing methods mentioned above may not be very efficient. Moreover, these methods are based on weakly admissible structured low-rank formats, i.e. HODLR and HSS-matrices, that achieve optimal complexity under the assumption that all off-diagonal blocks have small numerical rank that does not grow with the problem size. While this might be the case in some problems involving simple 1D or 2D geometries, for 3D problems in general, the off-diagonal blocks often have large ranks so using weakly admissible formats leads to suboptimal complexity. Therefore, for such problems, strongly admissible formats like $\mathcal{H}$ and $\mathcal{H}^2$-matrices are often preferred since they can achieve optimal complexity even for general 3D problems.

In this paper, we study a bisection method called slicing-the-spectrum that computes the $k$-th smallest eigenvalue of a symmetric matrix~\cite{Parlett1998}. This method relies on Sylverster's inertia theorem that allows one to compute the number of eigenvalues of a matrix $A$ that are smaller than a value $\mu$ by evaluating the LDL factorization of the shifted matrix $A - \mu I$.
It has been studied with HODLR matrices in~\cite{Benner2012_Bisection_HODLR}, allowing the computation of the $k$-th eigenvalue in $\mathcal{O}(n \log_2^4(n) \log_2((b-a)/\epsilon_{ev}))$ operations, where $[a, b]$ is the bisection starting interval and $\epsilon_{ev}$ is the desired accuracy. Generalized LDL factorization of HSS-matrices has also been used to compute the inertia in $\mathcal{O}(n)$, which ultimately reduces the cost down to $\mathcal{O}(n \log_2((b-a)/\epsilon_{ev}))$ per eigenvalue~\cite{Xi2014_Bisection_HSS}. Although this is optimal in terms of computational cost, the application to general 3D problems is very limited due to the weakness of HSS that makes the inertia evaluation no longer $\mathcal{O}(n)$. Nevertheless, it has been reported in~\cite{Xi2014_Bisection_HSS} that even when the off-diagonal block does not have a small numerical rank, aggressively discarding most of its singular values during HSS compression leads to a sufficiently accurate inertia evaluation when only a few largest eigenvalues are desired. However, when other parts of the spectrum are desired, especially the ones that are far from the largest eigenvalue, a large portion of the off-diagonal singular values need to be kept, leading to a high compression rank that grows with the problem size.

Here we present a slicing-the-spectrum method that uses a generalized LDL factorization of $\mathcal{H}^2$-matrices to compute the $k$-th eigenvalue in $\mathcal{O}(n \log_2((b-a)/\epsilon_{ev}))$ operations. Due to the flexible structure of $\mathcal{H}^2$-matrices that can handle dense off-diagonal blocks, our method is more efficient and applicable to a wider range of problems than the existing HODLR and HSS-based eigenvalue solvers. In addition, our method can also be used to compute some or even all eigenvalues. In such cases, some of the computed inertia can be reused to reduce the initial interval size when computing the other eigenvalues. Moreover, the bisection of disjoint intervals can be done in parallel since they are independent of each other. Therefore if some eigenvalues are desired, our method is still competitive with the existing QR iterations and divide-and-conquer eigenvalue solvers. However, when only a few interior eigenvalues are needed, our method is better since it can compute one eigenvalue in almost linear time without computing the others.

To the best of our knowledge, there is only one existing report that also uses $\mathcal{H}^2$-matrix factorization for slicing-the-spectrum~\cite{Benner2013_Bisection_H2}. While our method bears similarities with it, the key difference is in the LDL factorization used to slice the spectrum: the method in~\cite{Benner2013_Bisection_H2} uses a pure Cholesky-based LDL factorization with $\mathcal{O}(n \log_2(n))$ arithmetic complexity, 
whereas our method uses a generalized LDL factorization based on $\mathcal{H}^2$-ULV factorization~\cite{Ma2019_UMV} with $\mathcal{O}(n)$ complexity, meaning that our method requires less computational cost.

The rest of this paper is organized as follows. In Section 2 we introduce the slicing-the-spectrum method. In Section 3 we discuss our proposed linear time generalized LDL factorization of $\mathcal{H}^2$ matrices. Then in Section 4, we explain slicing the spectrum of $\mathcal{H}^2$-matrices with generalized LDL factorization and its parallelization. Numerical experiments that show the accuracy and efficiency of our method are presented in Section 5. Section 6 concludes the paper. 
\section{Slicing the Spectrum}
\label{section:slicing-the-spectrum}
In this section, we briefly recall the slicing-the-spectrum method~\cite{Parlett1998}. The main idea is to use bisection to find the $k$-th smallest eigenvalue $\lambda_k$ that is contained within the search interval $[a, b]$. In each step, the method evaluates the function $\nu(A - \mu I)$ which corresponds to the number of negative eigenvalues of the shifted matrix $A - \mu I$. Sylvester's inertia law tells us that this is equal to the number of negative entries in the diagonal matrix $D$ coming from the LDL factorization of the shifted matrix, i.e.
\begin{equation}
    \nu(A - \mu I) = \cardinality{\{ \lambda \in \Lambda(A) | \lambda < \mu \}}
    = \cardinality{\{ j | D_{j,j} < 0 \}},
\end{equation}
where $A - \mu I = LDL^T$. This function is used to choose the part of the interval that contains $\lambda_k$ in each iteration until the bisection process is stopped when the interval size is smaller than a prescribed threshold $\epsilon_{ev} > 0$. Algorithm~\ref{alg:slicing-the-spectrum} summarizes these steps.
\begin{algorithm}[h]
  \DontPrintSemicolon
  \caption{Slicing the Spectrum}
  \label{alg:slicing-the-spectrum}
  \SetAlgoLined
  \KwIn{$A \in \mathbb{R}^{n \times n}$, $k \in \mathbb{Z}$ ($1 \leq k \leq n$), and $a,b \in \mathbb{R}$ such that $\lambda_k \in [a, b]$}
  \KwOut{$\tilde{\lambda}_k$}
  $a_k = a; \; b_k = b$ \;
  \While{$b_k - a_k \geq \epsilon_{ev}$} {
    $\mu = (a_k + b_k)/2$ \;
    $A - \mu I = LDL^T$ \;
    $\nu(A - \mu I) = \cardinality{\{ j | D_{j,j} < 0 \}}$ \;
    \lIf{$\nu(A - \mu I) \geq k$} {
        $b_k \leftarrow \mu$
    }
    \lElse {
        $a_k \leftarrow \mu$
    }
  }
  $\tilde{\lambda}_k = (a_k + b_k) / 2$
\end{algorithm}

At the end of Algorithm~\ref{alg:slicing-the-spectrum}, we have the final interval $[a_k,b_k]$ that contains the $k$-th eigenvalue with $b_k - a_k < \epsilon_{ev}$, along with the approximate eigenvalue $\tilde{\lambda}_k = (a_k + b_k) / 2$. Thus, we have the upper bound of the eigenvalue error
\begin{equation}
\label{eq:error-bound}
    \abs{\lambda_k - \tilde{\lambda}_k} < \frac{1}{2} \epsilon_{ev}.
\end{equation}
This algorithm requires $\mathcal{O}(\log_2((b-a)/\epsilon_{ev}))$ iterations to reduce the interval size down to $\epsilon_{ev}$, where at each iteration one LDL factorization is computed. This is prohibitive for a general dense matrix $A$ that requires $\mathcal{O}(n^3)$ operations for LDL factorization. However, for structured low-rank matrices that allow fast LDL factorization, this algorithm becomes attractive, as we will see in the following sections.
\section{Generalized LDL Factorization of $\mathcal{H}^2$-Matrices}
\label{section:h2-ldlt}
In this section, we explain our proposed linear time generalized LDL factorization of $\mathcal{H}^2$-matrices, which is based on the $\mathcal{H}^2$-ULV factorization discussed in~\cite{Ma2019_UMV}. We start by introducing the hierarchical block matrix notation that we use throughout this paper. Then we describe the generalized LDL factorization for the weakly admissible BLR\textsuperscript{2}-matrices, which can be seen as a simple, non-hierarchical version of HSS. Next, we introduce hierarchy to obtain generalized HSS-LDL factorization. Finally, we extend that to strong admissibility to obtain a generalized $\mathcal{H}^2$-LDL factorization.

\subsection{Notation}
\begin{figure}[h]
    \centering
    \includegraphics[width=0.6\linewidth]{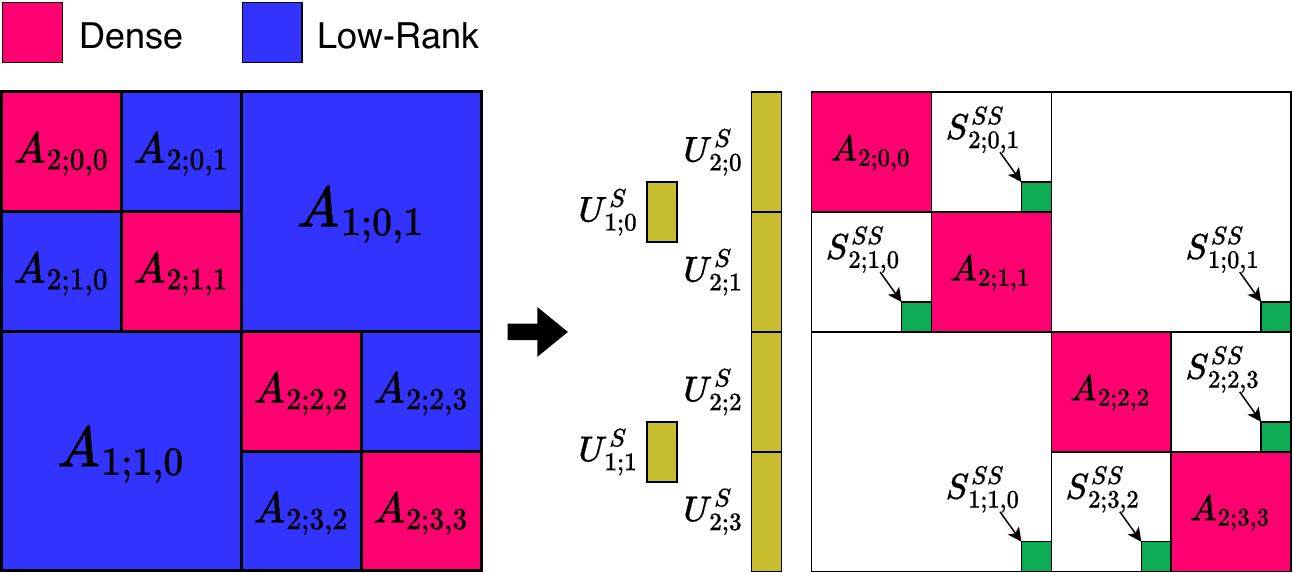}
    \caption{Index notation for hierarchical matrix block.}
    \label{fig:notation}
\end{figure}

Given an $n \times n$ symmetric matrix $A$, we index its hierarchically subdivided blocks as $A_{level;row,column}$. An approximated low-rank block is denoted as $\tilde{A}_{level;row,column}$. The shared column (row) basis is denoted as $U_{level;row}$. At the leaf level this denotes the actual shared basis, while at non-leaf levels it denotes the transfer matrix. A graphical illustration is shown in Fig.~\ref{fig:notation}.

For example, the block $\bluetext{A_{1;0,1}}$ can be approximated in hierarchical low-rank form as
\begin{equation}
\label{eq:notation-example}
    \bluetext{\tilde{A}_{1;0,1}} =
    \yellowtext{
        \begin{bmatrix}
            U_{2;0}^S & 0 \\
            0 & U_{2;1}^S
        \end{bmatrix}
        U_{1;0}^S
    }
    \greentext{S_{1;0,1}^{SS}}
    \yellowtext{
        \tp{U_{1;1}^S}
        \begin{bmatrix}
            \tp{U_{2;2}^S} & 0 \\
            0 & \tp{U_{2;3}^S}
        \end{bmatrix}
    },
\end{equation}
where $\greentext{S_{1;0,1}^{SS}}$ is the small \greentext{skeleton matrix} shown inside the large off-diagonal block on the top-right of Fig~\ref{fig:notation}. The shared column basis for a given row is obtained by performing the rank-revealing QR factorization, e.g.
\begin{align}
    \yellowtext{
        \begin{bmatrix}
            U_{2;0}^{S} & U_{2;0}^{R}
        \end{bmatrix}
    }, R &= QR(\yellowtext{A_{2;0,+}}), \label{eq:shared-bases} \\
    \yellowtext{
        \begin{bmatrix}
            U_{1;0}^S & U_{1;0}^R
        \end{bmatrix}
    }, R &= QR\left(
        \yellowtext{
            \begin{bmatrix}
                U_{2;0}^S & 0 \\
                0 & U_{2;1}^S
            \end{bmatrix}^{\top}
            A_{1;0,+}
        }
    \right), \label{eq:transfer-matrices}
\end{align}
where $A_{level;row,+}$ denotes the concatenation of all low-rank blocks in the entire row at a particular level. Note that other rank revealing factorization like Singular Value Decomposition or Interpolative Decomposition~\cite{HalkoMartinsson2011_LRA} can also be used here. The $S$ and $R$ superscripts denotes the skeleton and redundant part of the basis, respectively. The skeleton part corresponds to the approximate column basis, e.g.
\begin{equation}
    \norm{\yellowtext{U_{2;0}^S \tp{U_{2;0}^S} A_{2;0,+} - A_{2;0,+}}} \leq \epsilon_{\mathcal{H}^2}
\end{equation}
for a prescribed error tolerance $\epsilon_{\mathcal{H}^2}$, while the redundant is part of the basis that is usually discarded during the approximation.
The \yellowtext{shared bases} along with the \greentext{skeleton matrices} for all low-rank blocks can be computed in $\mathcal{O}(n)$ using $\mathcal{H}^2$ construction techniques discussed in~\cite{Borm2005_HCA,Borm2003_Interpolation,Cai2022_HiDR}. Note that to perform the generalized LDL factorization, we assume the shared bases to have orthonormal columns. Orthogonalization of the shared bases can also be done in $\mathcal{O}(n)$ if necessary.

\subsection{Generalized BLR\textsuperscript{2}-LDL Factorization}
\begin{figure}[h]
    \centering
    \includegraphics[width=0.6\linewidth]{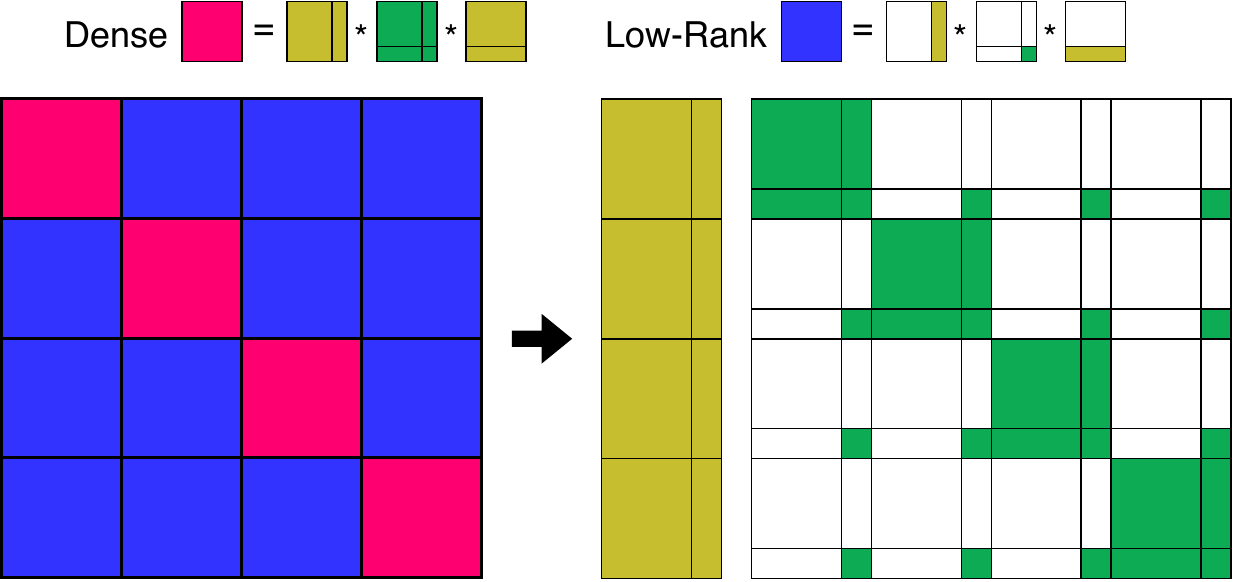}
    \caption{Introducing zeros to low-rank blocks.}
    \label{fig:blr2-usv}
\end{figure}

Given a weakly admissible BLR\textsuperscript{2}-matrix, the generalized LDL factorization is described as follows. First, we form the \greentext{skeleton matrices} for both dense and low-rank blocks by incorporating the redundant part of the shared bases. For each dense diagonal block $\redtext{A_{2;i,i}}$, the skeleton matrix is obtained from
\begin{equation}
\label{eq:dense_skeleton}
    \greentext{
        \begin{bmatrix}
            S_{2;i,i}^{RR} & S_{2;i,i}^{RS} \\
            S_{2;i,i}^{SR} & S_{2;i,i}^{SS}
        \end{bmatrix}
    } =
    \yellowtext{
        \begin{bmatrix}
            U_{2;i}^R & U_{2;i}^S
        \end{bmatrix}^{\top}
    }
    \redtext{A_{2;i,i}} \;
    \yellowtext{
        \begin{bmatrix}
            U_{2;i}^R & U_{2;i}^S
        \end{bmatrix}
    }.
\end{equation}
For each low-rank block $\bluetext{\tilde{A}_{2;i,j}}$ ($i \neq j$), the skeleton matrix is obtained by
\begin{align}
\label{eq:lr_skeleton}
    \greentext{
        \begin{bmatrix}
            S_{2;i,j}^{RR} & S_{2;i,j}^{RS} \\
            S_{2;i,j}^{SR} & S_{2;i,j}^{SS}
        \end{bmatrix}
    } &=
    \yellowtext{
        \begin{bmatrix}
            U_{2;i}^R & U_{2;i}^S
        \end{bmatrix}^{\top}
    }
    \bluetext{\tilde{A}_{2;i,j}} \;
    \yellowtext{
        \begin{bmatrix}
            U_{2;j}^R & U_{2;j}^S
        \end{bmatrix}
    } \nonumber \\
    & =
    \greentext{
        \begin{bmatrix}
            0 & 0 \\
            0 & \yellowtext{\tp[-6]{U_{2;i}^S}}
            \bluetext{A_{2;i,j}} \yellowtext{U_{2;j}^S}
        \end{bmatrix}
    }.
\end{align}
Now the basis $\yellowtext{\begin{bmatrix} U_{2;i}^R & U_{2;i}^S \end{bmatrix}}$ and $\yellowtext{\begin{bmatrix} U_{2;j}^R & U_{2;j}^S \end{bmatrix}^{\top}}$ are shared among both dense and low-rank blocks in the same row and column, respectively. Consequently, this introduces zeros to the low-rank blocks, as written in Eq.(\ref{eq:lr_skeleton}) and graphically shown in Fig.~\ref{fig:blr2-usv}.

\begin{figure*}[h]
    \centering
    \includegraphics[width=0.8\linewidth]{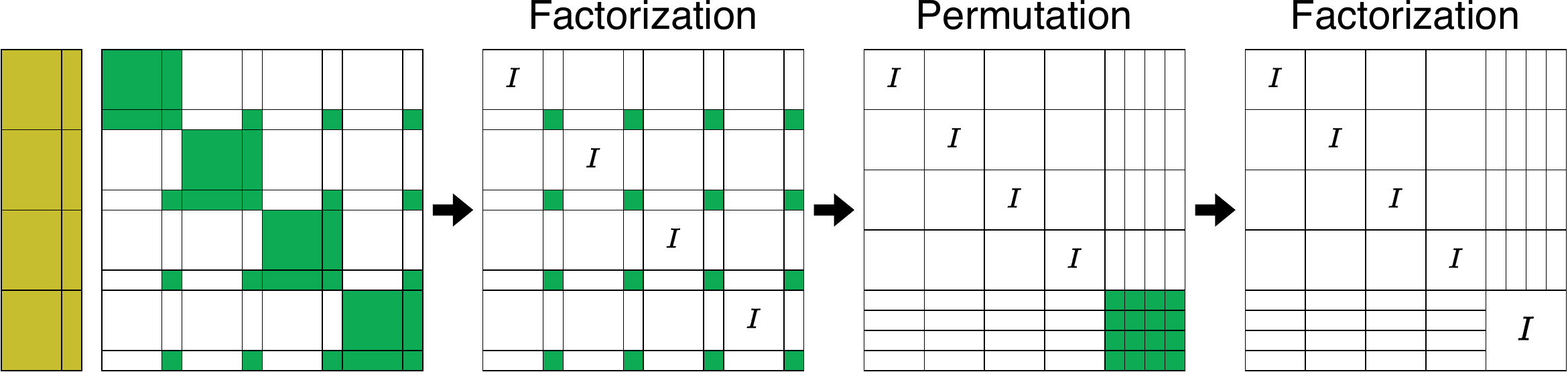}
    \caption{Flow of Generalized BLR\textsuperscript{2}-LDL Factorization.}
    \label{fig:blr2-ldl}
\end{figure*}
The next step is to perform factorization of the skeleton matrices, as described in Fig~\ref{fig:blr2-ldl}. First, we compute the LDL factorization of the redundant part of the diagonal blocks to obtain
\begin{equation}
\label{eq:blr2-diag-ldl}
     L_{2;i,i}^{RR}, D_{2;i,i}^{RR} = LDL^T(S_{2;i,i}^{RR}).
\end{equation}
Then we eliminate the $S_{2;i,i}^{SR}$ blocks
\begin{equation}
\label{eq:blr2-diag-elimination}
    L_{2;i,i}^{SR} = S_{2;i,i}^{SR} 
    \left(L_{2;i,i}^{RR} \right)^{-T}
    \left(D_{2;i,i}^{RR} \right)^{-1},
\end{equation}
followed by computing the Schur's complements
\begin{equation}
\label{eq:blr2-diag-schur}
    S_{2;i,i}^{SS} = S_{2;i,i}^{SS} - L_{2;i,i}^{SR} D_{2;i,i}^{RR}
    \tp[-6]{L_{2;i,i}^{SR}}.
\end{equation}
The processes described in Eqs.(\ref{eq:blr2-diag-ldl})-(\ref{eq:blr2-diag-schur}) for different $i$ can be done concurrently since there are no dependencies among them, as shown in the first factorization phase in Fig.~\ref{fig:blr2-ldl}. Once the $S_{2;i,i}^{RR}$ part of the diagonal blocks are factorized, the remaining $S_{2;i,j}^{SS}$ blocks to be factorized are permuted to the bottom-right corner, as shown in the permutation phase in Fig.~\ref{fig:blr2-ldl}. Finally, the entire remaining block is factorized as a single dense matrix:
\begin{equation}
\label{eq:blr2-root-ldl}
    L_{2;0:3,0:3}^{SS}, D_{2;0:3,0:3}^{SS} =
    LDL^T\left(
    \begin{bmatrix}
        S_{2;0,0}^{SS} & S_{2;0,1}^{SS} & S_{2;0,2}^{SS} & S_{2;0,3}^{SS} \\
        S_{2;1,0}^{SS} & S_{2;1,1}^{SS} & S_{2;1,2}^{SS} & S_{2;1,3}^{SS} \\
        S_{2;2,0}^{SS} & S_{2;2,1}^{SS} & S_{2;2,2}^{SS} & S_{2;2,3}^{SS} \\
        S_{2;3,0}^{SS} & S_{2;3,1}^{SS} & S_{2;3,2}^{SS} & S_{2;3,3}^{SS}
    \end{bmatrix}
    \right)
\end{equation}

\subsection{Generalized HSS-LDL Factorization}
\begin{figure*}[h]
    \centering
    \includegraphics[width=0.8\linewidth]{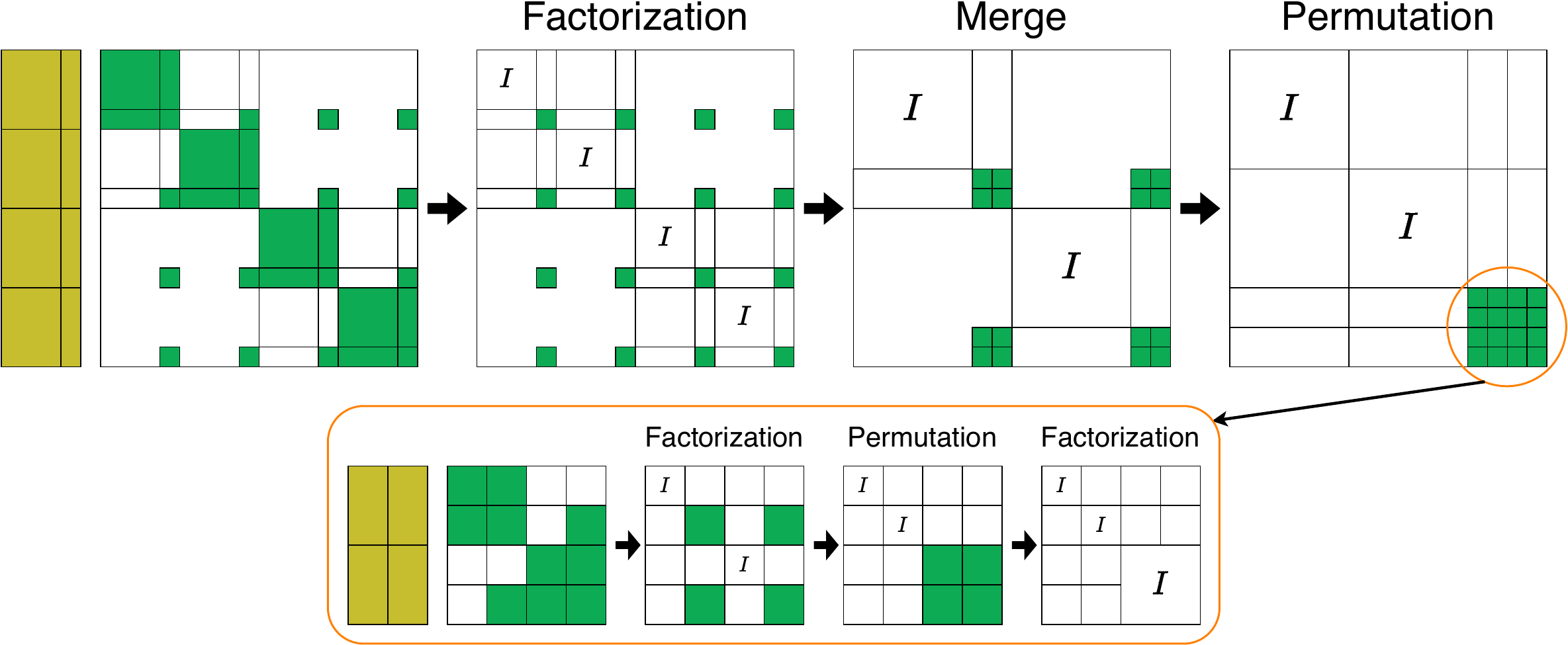}
    \caption{Flow of Generalized HSS-LDL Factorization.}
    \label{fig:hss-ldl}
\end{figure*}

Here we introduce hierarchy and describe the generalized LDL factorization of HSS matrices, which is a multi-level version of the weakly admissible BLR\textsuperscript{2} matrices. In this case, the factorization proceeds in a bottom-up fashion, starting from the leaf level and going up to the upper levels. In order to illustrate this, we use an example of a 2-level HSS matrix as shown in Fig.~\ref{fig:hss-ldl}. The same recursive steps can be applied to general HSS matrices of arbitrary size with any number of subdivision levels.

Fig.~\ref{fig:hss-ldl} shows that the flow at each level is identical to that of the BLR\textsuperscript{2} version. Starting at the leaf level ($level=2$), we perform a similar process described in Eqs.(\ref{eq:dense_skeleton})-(\ref{eq:blr2-diag-schur}). Then instead of factorizing the entire remaining blocks as a single dense matrix shown in Eq.(\ref{eq:blr2-root-ldl}), we treat them as another BLR\textsuperscript{2}-matrix (or HSS-matrix whose level is reduced by 1) and recursively apply a similar procedure. At this level, the remaining blocks are partitioned as
\begin{equation}
    \begin{bmatrix}
        \redtext{\hat{S}_{1;0,0}} & \bluetext{\hat{S}_{1;0,1}} \\
        \bluetext{\hat{S}_{1;1,0}} & \redtext{\hat{S}_{1;1,1}}
    \end{bmatrix} =
    \begin{bmatrix}
        S_{2;0,0}^{SS} & S_{2;0,1}^{SS} & S_{2;0,2}^{SS} & S_{2;0,3}^{SS} \\
        S_{2;1,0}^{SS} & S_{2;1,1}^{SS} & S_{2;1,2}^{SS} & S_{2;1,3}^{SS} \\
        S_{2;2,0}^{SS} & S_{2;2,1}^{SS} & S_{2;2,2}^{SS} & S_{2;2,3}^{SS} \\
        S_{2;3,0}^{SS} & S_{2;3,1}^{SS} & S_{2;3,2}^{SS} & S_{2;3,3}^{SS}
    \end{bmatrix},
\end{equation}
as shown in the merge and permutation steps at the top part of Fig.~\ref{fig:hss-ldl}. Also, now the shared bases $\yellowtext{U_{1;i}}$ are the transfer matrices along with their redundant parts as shown in Eq.(\ref{eq:transfer-matrices}). Next, we use these shared bases to form the \greentext{skeleton matrices} of both dense and low-rank blocks, similar to that of Eqs.(\ref{eq:dense_skeleton})-(\ref{eq:lr_skeleton}). For each dense diagonal block $\redtext{\hat{S}_{1;i,i}}$ we obtain
\begin{equation}
\label{eq:hss_dense_skeleton}
    \greentext{
        \begin{bmatrix}
            S_{1;i,i}^{RR} & S_{1;i,i}^{RS} \\
            S_{1;i,i}^{SR} & S_{1;i,i}^{SS}
        \end{bmatrix}
    } =
    \yellowtext{
        \begin{bmatrix}
            U_{1;i}^R & U_{1;i}^S
        \end{bmatrix}^{\top}
    }
    \redtext{\hat{S}_{1;i,i}}
    \yellowtext{
        \begin{bmatrix}
            U_{1;i}^R & U_{1;i}^S
        \end{bmatrix}
    },
\end{equation}
whereas for off-diagonal low-rank block $\bluetext{\hat{S}_{1;i,j}}$ ($i \neq j$) we compute
\begin{align}
\label{eq:hss_lr_skeleton}
    \greentext{
        \begin{bmatrix}
            S_{1;i,j}^{RR} & S_{1;i,j}^{RS} \\
            S_{1;i,j}^{SR} & S_{1;i,j}^{SS}
        \end{bmatrix}
    } &=
    \yellowtext{
        \begin{bmatrix}
            U_{1;i}^R & U_{1;i}^S
        \end{bmatrix}^{\top}
    }
    \bluetext{\hat{S}_{1;i,j}}
    \yellowtext{
        \begin{bmatrix}
            U_{1;j}^R & U_{1;j}^S
        \end{bmatrix}
    } \nonumber \\
    & =
    \greentext{
        \begin{bmatrix}
            0 & 0 \\
            0 & \yellowtext{\tp[-6]{U_{1;i}^S}}
            \bluetext{\hat{S}_{1;i,j}} \yellowtext{U_{1;j}^S}
        \end{bmatrix}
    }.
\end{align}
Finally, once we have the skeleton matrices at this level, we apply the same process shown in Eqs.(\ref{eq:blr2-diag-ldl})-(\ref{eq:blr2-root-ldl}), as shown in the bottom part of Fig.~\ref{fig:hss-ldl}.

\subsection{Generalized $\mathcal{H}^2$-LDL Factorization}
\begin{figure*}[h]
    \centering
    \includegraphics[width=0.95\linewidth]{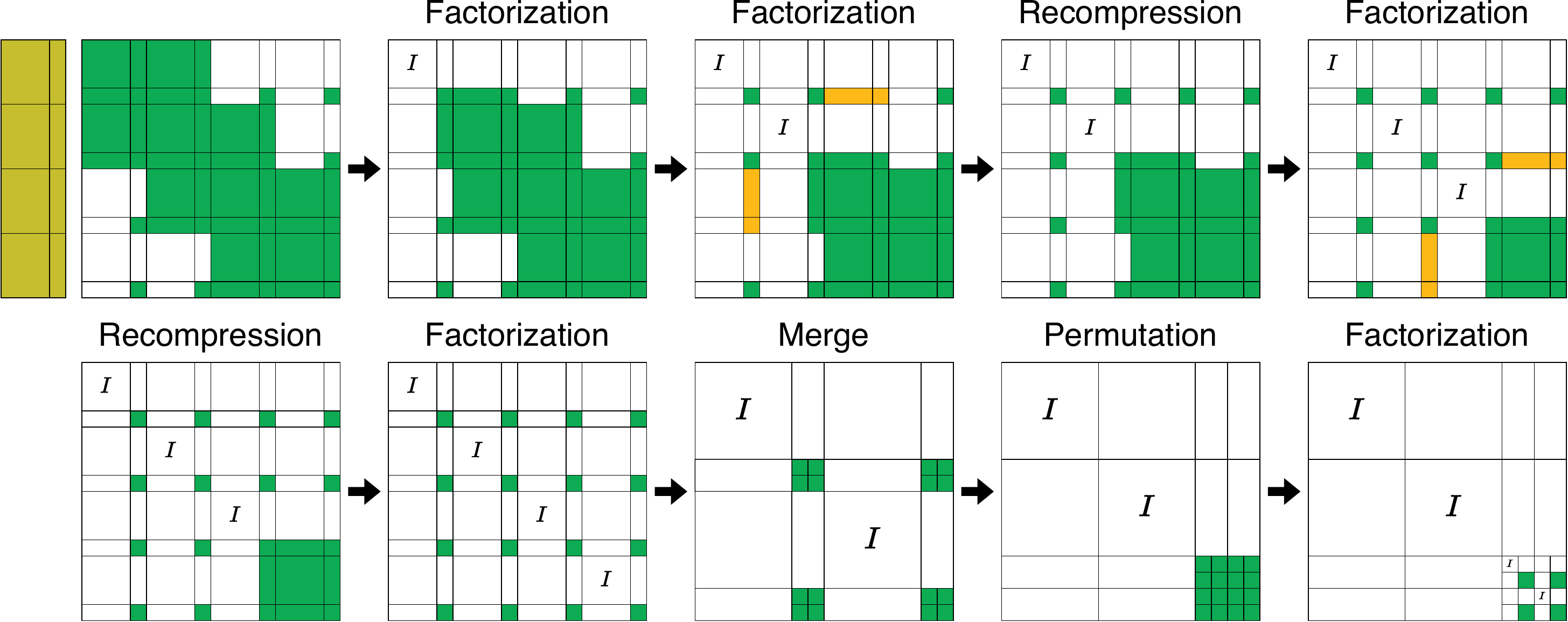}
    \caption{Flow of Generalized $\mathcal{H}^2$-LDL Factorization.}
    \label{fig:h2-ldl}
\end{figure*}

In this section, we extend the generalized LDL factorization for HSS to $\mathcal{H}^2$ matrices by introducing dense off-diagonal blocks. We take as an example a 2-level $\mathcal{H}^2$ matrix $\tilde{A}$ with dense off-diagonal blocks as shown in the leftmost part of Fig.~\ref{fig:h2-ldl}. While this example has a rather simple block tri-diagonal structure, the same steps can be applied to general $\mathcal{H}^2$ matrices with arbitrary subdivision levels and any pattern of dense off-diagonal blocks.

The flow at each level begins by forming the skeleton matrices using both the redundant and skeleton part of the shared bases, which is identical to the BLR\textsuperscript{2} version described in Eqs.(\ref{eq:dense_skeleton})-(\ref{eq:lr_skeleton}). Note that dense off-diagonal blocks are never used for the construction of shared bases. Also, the skeleton matrix of the dense off-diagonal block takes a similar form as the dense diagonal block, meaning that Eq.(\ref{eq:dense_skeleton}) also applies to the dense off-diagonal blocks.

The next step is to factorize the redundant part of the skeleton matrices. As can be seen in Fig.~\ref{fig:h2-ldl}, dense off-diagonal blocks lead to the computation of Schur's complements in the trailing submatrices. When it fills into a low-rank block, it creates a \orangetext{fill-in} block that may destroy the sparsity of the matrix, shown as the orange block in Fig.~\ref{fig:h2-ldl}. In order to mitigate this, the fill-ins in an entire row/column $i$ are incorporated into the shared basis just before the $i$-th diagonal block is factorized. This is the recompression step where the orange blocks are merged into the shared basis and disappear into the low-rank block. For example in the third matrix from the left on the top part of Fig.~\ref{fig:h2-ldl}, there are two \orangetext{fill-in} blocks within $\greentext{S_{2;0,2}}$ and $\greentext{S_{2;2,0}}$. So before factorizing $\greentext{S_{2;2,2}^{RR}}$, a recompression step is performed that involves an update to the shared basis
\begin{displaymath}
    \yellowtext{
        \begin{bmatrix}
            U_{2;2}^S & U_{2;2}^R
        \end{bmatrix}
    }, R \leftarrow QR \left(
        \begin{bmatrix}
            \yellowtext{U_{2;2}^S} \greentext{S_{2;2,0}^{SS}} &
            \yellowtext{
                \begin{bmatrix}
                    U_{2;2}^R & U_{2;2}^S
                \end{bmatrix}
            }
            \begin{bmatrix}
                \orangetext{F_{2;2,0}^{RS}} \\
                \orangetext{F_{2;2,0}^{SS}}
            \end{bmatrix}
        \end{bmatrix}
    \right),
\end{displaymath}
where $\orangetext{F_{2;2,0}}$ is the corresponding fill-in block within $\greentext{S_{2;2,0}}$. This incorporates the fill-ins to the shared column basis so that they can be merged into $\greentext{S_{2;0,2}^{SS}}$ and $\greentext{S_{2;2,0}^{SS}}$, thus keeping the sparsity of the \greentext{skeleton} matrices intact. Note that in this case there is only one low-rank block in the particular row/column. However, in general it is necessary to concatenate the products for all low-rank blocks in the row/column that is being recompressed, which for example can be done efficiently using the technique described in~\cite{Ma2019_UMV}.

Once the redundant part of the skeleton matrices are factorized, we merge and permute the remaining blocks then recursively apply the same procedure until the top level is reached and we are left with a single dense matrix, similarly to the HSS version. The overall cost of the algorithm is $\mathcal{O}(n)$ as long as the numerical ranks of the off-diagonal low-rank blocks grow independently of the matrix size (see the discussion in~\cite{Ma2019_UMV}). At the end of the whole process, we are left with a factorization in the form of
\begin{equation}
\label{eq:h2-generalized-ldl}
    \tilde{A} \approx \mathcal{LDL}^T,
\end{equation}
where $\mathcal{L}$ is a product of lower triangular and orthogonal matrices and $\mathcal{D}$ is a diagonal matrix.

\section{Slicing the Spectrum with Generalized $\mathcal{H}^2$-LDL Factorization}
\label{section:h2-slicing-the-spectrum}
Sylvester's theorem tells us that the inertia of a matrix is invariant under congruence transformations. 
Therefore, based on Eq.(\ref{eq:h2-generalized-ldl}) the inertia of $\tilde{A}$ can also be obtained by looking at the diagonal entries of $\mathcal{D}$. So we can use Algorithm~\ref{alg:slicing-the-spectrum} to compute the $k$-th smallest eigenvalue of $\tilde{A}$ by replacing the LDL factorization of the shifted matrix in line 4 with the generalized LDL factorization shown in Eq.(\ref{eq:h2-generalized-ldl}). This gives us a cost of $\mathcal{O}(n)$ per iteration to compute the inertia of the shifted matrix $\tilde{A} - \mu I$, leading to the total cost of
\begin{equation}
    \mathcal{O}(n \log_2((b-a)/\epsilon_{ev}))
\end{equation}
operations to compute the $k$-th eigenvalue with an accuracy of $\epsilon_{ev}$. Moreover, this process only requires additional storage of one $\mathcal{H}^2$ matrix to store the shifted matrix at each iteration, which implies an overall storage requirement of $\mathcal{O}(n)$.

\subsection{Accuracy}
It is important to note that our proposed generalized LDL factorization shown in Eq.(\ref{eq:h2-generalized-ldl}) is approximative. Therefore, we must make sure that the factorization is sufficiently accurate to produce the correct inertia evaluation of the shifted matrices. According to~\cite{Parlett1998}, the generalized LDL factorization is sufficiently accurate for inertia evaluation if it satisfies
\begin{equation}
    \norm{\left(\tilde{A}-\mu I\right) - \mathcal{L_{\mu} D_{\mu} L_{\mu}}^T} \leq \underset{j}{\text{min}} \abs{\lambda_j(\tilde{A})-\mu},
\end{equation}
where $\lambda_j(\tilde{A})$ is the $j$-th eigenvalue of $\tilde{A}$ and $\tilde{A}-\mu I \approx \mathcal{L_{\mu} D_{\mu} L_{\mu}}^T$ is the factorization of the shifted matrix. Therefore, we need an error bound of the form
\begin{equation}
    \norm{\tilde{A}- \mathcal{LDL}^T} \leq \epsilon.
\end{equation}
In addition, we also need this bound for all shifted matrices $\tilde{A}-\mu I$. To the best of our knowledge, such a bound is not available in the current literature on hierarchical matrices, especially for an $\mathcal{H}^2$ matrix $\tilde{A}$. Although it has been mentioned in~\cite{Ma2019_UMV} that the accuracy of the similar $\mathcal{H}^2$-ULV factorization is directly controlled by the truncation rank used in the recompression step, there has been no discussion about the theoretical error bound. Thus at the moment, we do not have theoretical proof that the proposed generalized $\mathcal{H}^2$-LDL factorization is accurate for slicing the spectrum in general. Fortunately, many numerical experiments have shown that $\mathcal{H}^2$-ULV factorization produces satisfying accuracy~\cite{Ma2018_UMV,Ma2019_UMV,Ma2022_ULV}. Thus, we provide evidence that our method is sufficiently accurate through the numerical experiments shown in Section~\ref{section:results}.

\subsection{Parallelization}
\label{subsec:parallelization}
If only one eigenvalue is desired, one possible location to exploit parallelism is the generalized $\mathcal{H}^2$-LDL factorization. However, in our current method, the dependencies between the fill-in recompression and diagonal factorization limit the parallelization to the triangular solves and Schur complements involving dense blocks at each level. Since the number of such dense blocks usually grows somewhat constantly independent of the matrix size, the degree of parallelism is limited.

On the other hand, if more than one eigenvalue is desired, we can parallelize the bisection part to compute them efficiently. There are two key ideas for this:
\begin{itemize}
    \item Some of the computed inertia can be reused for the computation of other eigenvalues. For example, in the beginning, the same starting interval $[a, b]$ is assumed for the computation of all target eigenvalues. Suppose after the first inertia evaluation we obtain $\nu(A - \mu I)=v$. Then we know that the interval $[a, \mu)$ contains the first $v$ eigenvalues and $[\mu, b]$ contains the remaining $n-v$ eigenvalues. Thus, for the computation of the $k$-th eigenvalue where $k \leq v$, we can narrow the starting interval to $[a, \mu)$. The same applies for $k > v$ where the starting interval of $[\mu, b]$ can be used.
    \item The bisection of disjoint intervals are independent of each other so they can be done in parallel. For this, each process needs to hold an instance of the matrix along with the interval of interest. Although this is quite prohibitive for dense matrices that require $\mathcal{O}(n^2)$ storage, for $\mathcal{H}^2$-matrices that require only $\mathcal{O}(n)$ storage this is still feasible even for relatively large matrices.
\end{itemize}
Thus, given the starting interval $[a, b]$ and a range of eigenvalue indices $[k_0, k_1]$ to be computed, we consider two approaches to compute the eigenvalues in parallel on distributed memory systems: one that distributes the eigenvalue indices equally to every process, and another that uses a master-worker model. In both approaches, we assume that every process holds an instance of the $\mathcal{H}^2$ matrix whose eigenvalues will be computed.

In the first parallel algorithm, we distribute the eigenvalue index range $[k_0, k_1]$ equally to $P$ processes. Every process will first compute the smallest and largest eigenvalue from the set it is responsible for. After that, it computes the remaining inner eigenvalues using the smallest and largest eigenvalues as the starting interval. This method does not involve any communication during the eigenvalue computation since every process knows the set of indices it has to compute. However, this method may suffer from load imbalance when the eigenvalues are not uniformly distributed along the spectrum since some processes may be responsible for larger intervals than others even though the number of eigenvalues to be computed are the same. Hence, this method works best if the number of processes is close to the total number of eigenvalues to be computed. That way, every process would be responsible for a smaller part of the spectrum, thus reducing the overall load imbalance.

For the second algorithm, we adopt a master-worker model similar to the one discussed in~\cite{Mach2012_PhDThesis}. The master distributes the works by giving each idle worker an interval and a range of eigenvalue indices. The worker computes all eigenvalues within the given interval if there are no more than $2m$ eigenvalues to be computed. Otherwise, the worker will only compute the inner $m$ eigenvalues. The worker first computes the smallest and the largest eigenvalues within the set and notify the master that it will not compute the eigenvalues within the other two intervals, i.e. $[a, \tilde{\lambda}_{k}]$ and $[\tilde{\lambda}_{k+m-1}, b]$ where $[a, b]$ is the interval received by the worker and $\tilde{\lambda}_{k}$ and $\tilde{\lambda}_{k+m-1}$ are the smallest and largest eigenvalues from the chosen set of $m$ inner eigenvalues, respectively. After that, the worker continues to compute the rest of the eigenvalues in the set. The master manages a list of all free intervals and distributes them immediately to idle workers. The value $m$ is set to be small enough to balance the load among workers. This approach works well when computing a large number of eigenvalues, as reported in~\cite{Mach2012_PhDThesis,Benner2013_Bisection_H2}.
\section{Results}
\label{section:results}
In this section, we use several examples to demonstrate the performance and accuracy of our proposed method. Implementations were written in C++ where floating point calculations were performed in double precision. Single-threaded BLAS and LAPACK routines from Intel MKL were used for inner kernels involving dense matrices. Distributed memory parallelism was performed using MPI. Experiments were conducted on a system described in Table \ref{table:machine-specs}.
\begin{table}[hbtp]
    \centering
    \caption{Details of compute node used for experiments}
    \label{table:machine-specs}
    \begin{tabular}{ll}
        \toprule
        CPU & 2x AMD Ryzen Threadripper 3960X \\
        Clock speed & 3.8 GHz \\
        \# cores & 2 x 24 = 48 \\
        Memory & 120 GB \\
        Compiler suite & GCC 9.4 \\
        BLAS/LAPACK \& MPI & Intel 2022.1.0 \\
        \bottomrule
    \end{tabular}
\end{table}

The following methods have been compared:
\begin{itemize}
    \item $\mathbf{\mathcal{H}^2}$\textbf{-Bisection}: Slicing the spectrum with our generalized $\mathcal{H}^2$-LDL factorization described in Section~\ref{section:h2-slicing-the-spectrum}.
    \item \textbf{HSS-Bisection}: Slicing the spectrum with generalized HSS-LDL factorization~\cite{Xi2014_Bisection_HSS}.
    \item \textbf{LAPACK {\ttfamily{dsyevx}}}: Dense symmetric eigenvalue solver that computes selected eigenvalue(s) with controllable accuracy $\epsilon_{ev}$~\cite{Lapack1999}.
    \item \textbf{ScaLAPACK \ttfamily{pdsyev}}: Distributed-memory parallel dense symmetric eigenvalue solver that computes all eigenvalues in double precision accuracy~\cite{Scalapack1997}.
    \item \textbf{ELPA}: Distributed memory parallel dense eigenvalue solver that computes all eigenvalues using a 2-stage solver~\cite{Marek2014_ELPA}.
\end{itemize}
In order to ensure that the bisection methods produce accurate inertia evaluation, we used the low-rank compression threshold that was smaller than the bisection threshold, e.g. $\epsilon_{\mathcal{H}^2} = 10^{-2} \epsilon_{ev}$. The eigenvalue error was obtained by comparing the approximated eigenvalue with the result from LAPACK dense symmetric eigenvalue solver {\ttfamily{dsyev}}.
\subsection{Computing Some or All Eigenvalues of Synthetic Matrices}
Here we used the Laplace kernel to generate the rank-structured matrix
\begin{displaymath}
    A_{i,j} = \frac{1}{|x_i - x_j| + 10^{-3}},
\end{displaymath}
where $|x_i - x_j|$ denotes the Euclidean distance between two points $x_i$ and $x_j$. The points $x_i$ were uniformly distributed on the circumference of a unit circle.

\begin{figure}[hbtp]
    \centering
    \includegraphics[width=0.5\textwidth]{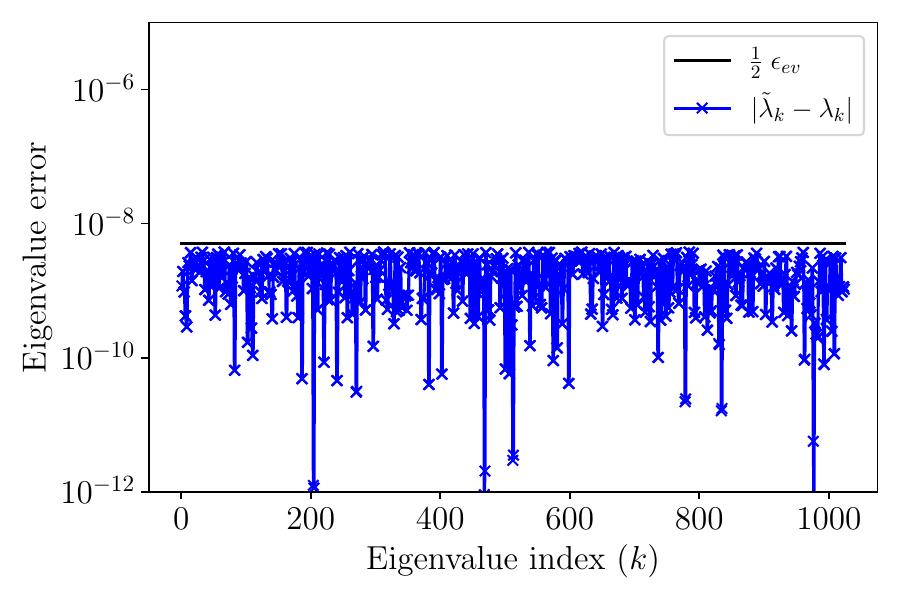}
    \caption{Eigenvalue errors for a $1024 \times 1024$ synthetic Laplace matrix with $\epsilon_{ev}=10^{-8}$.}
    \label{fig:laplace-acc-ev8}
\end{figure}
\begin{figure}[hbtp]
    \centering
    \includegraphics[width=0.5\textwidth]{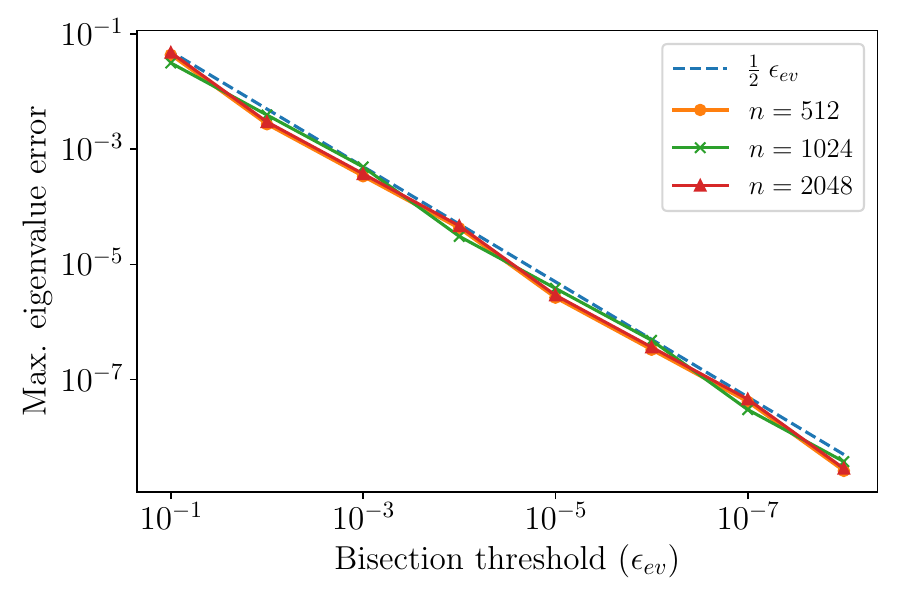}
    \caption{Maximum eigenvalue errors when computing all eigenvalues of synthetic Laplace matrices.}
    \label{fig:laplace-acc-all}
\end{figure}
We first show that our method produces accurate eigenvalues according to the bisection threshold. If every inertia evaluation is sufficiently accurate, the binary search will choose the correct interval that contains the target eigenvalue. Otherwise, it will choose the wrong interval that does not contain the target eigenvalue, leading to an error larger than the theoretical bound in Eq.(\ref{eq:error-bound}). Fig.~\ref{fig:laplace-acc-ev8} shows that all of the eigenvalue errors were below the theoretical error bound. Note that the eigenvalue error corresponds to the distance of the actual eigenvalue to the midpoint of our bisection interval in the last iteration. So different choices of starting interval may give slightly different errors, and in some cases, might result in an eigenvalue error that is much smaller than the desired accuracy, as shown in Fig.~\ref{fig:laplace-acc-ev8} where some of the eigenvalue errors were smaller than $10^{-10}$ even when $\epsilon_{ev}=10^{-8}$. Fig.~\ref{fig:laplace-acc-all} also shows that our method consistently produced accurate eigenvalues as we decreased the bisection threshold. We did not observe high rank growth that is proportional to the matrix size when using shifts near the actual eigenvalues as reported in~\cite{Benner2012_Bisection_HODLR} for $\mathcal{H}$-LDL factorization, but we do not have a theoretical bound for the ranks.

\begin{figure}[hbtp]
    \centering
    \begin{subfigure}[b]{0.49\textwidth}
        \centering
        \includegraphics[width=\linewidth]{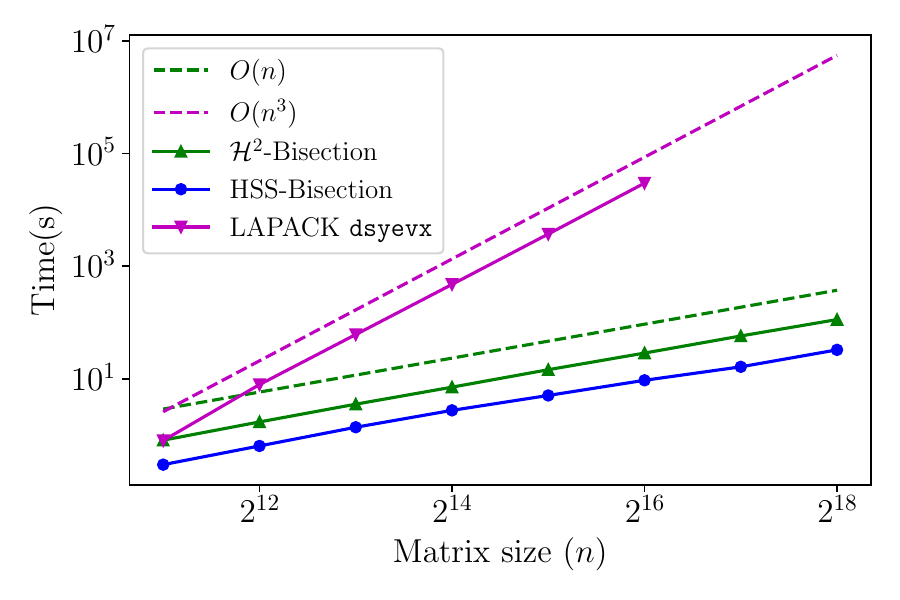}
        \caption{}
        \label{fig:laplace-eig-time}
    \end{subfigure}
    \hfill
    \begin{subfigure}[b]{0.49\textwidth}
        \centering
        \includegraphics[width=\linewidth]{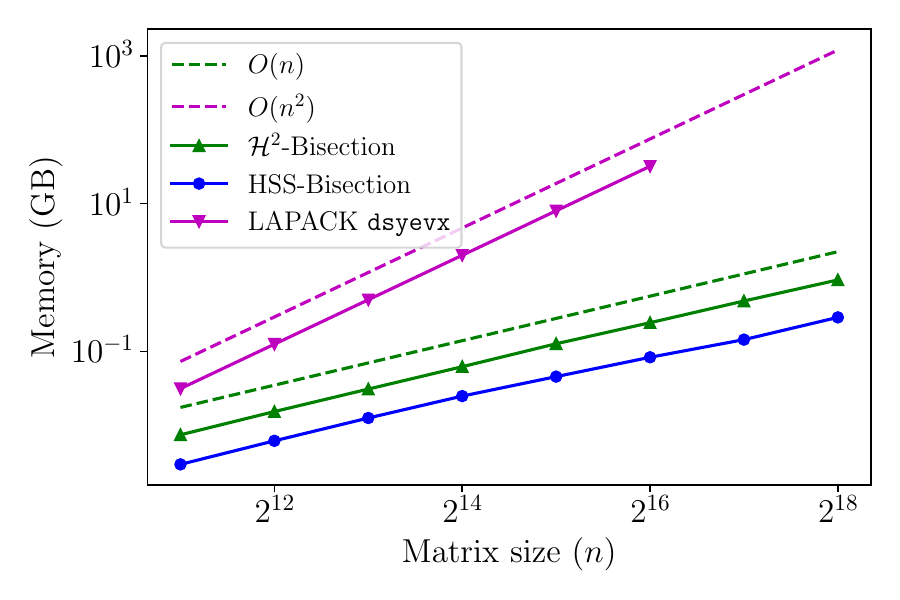}
        \caption{}
        \label{fig:laplace-eig-mem}
    \end{subfigure}
    \caption{Computing the $k$-th eigenvalue of synthetic Laplace matrices with $\epsilon_{ev}=10^{-5}$ using a single core of the compute node. (a) Computation time. (b) Memory consumption.}
    \label{fig:laplace-eig-res}
\end{figure}
Next, we show the performance of our method to compute the $k$-th eigenvalue where $k = n/2$ in Fig.~\ref{fig:laplace-eig-res}. A constant starting interval of $[0, 2048]$ was used for the bisection methods, which we found still contain the target (median) eigenvalue even if we increased the matrix size. Fig.~\ref{fig:laplace-eig-res} shows that with constant accuracy and size of starting interval, our $\mathcal{H}^2$ bisection method scaled linearly with respect to the matrix size in terms of both computation time and memory consumption, which is in accordance with our estimate in Section~\ref{section:h2-slicing-the-spectrum}. The HSS bisection method also showed a similar linear scaling since the test matrices originate from simple 2D geometries where HSS could still maintain a constant rank as we increased the problem size. We observed a maximum compression rank of 17 for the $\mathcal{H}^2$-Bisection and 50 for the HSS-Bisection. With this small difference in the compression rank, our method was slightly slower than the HSS one due to the hidden constant number of operations coming from the fill-in recompression steps that did not occur in generalized HSS-LDL factorization. Moreover, both bisection methods outperformed LAPACK {\ttfamily{dsyevx}} in terms of computation time and memory consumption. For the matrix of order 65,536, our bisection method was already about 3 orders of magnitude faster than LAPACK {\ttfamily{dsyevx}}. For the large matrix of order $n =  262,144$, our method required only about 0.28 GB of memory, whereas LAPACK {\ttfamily{dsyevx}} required about 512 GB, which exceeded the memory capacity per node of the system that we were using.

\begin{figure}[hbtp]
    \centering
    \begin{subfigure}[b]{0.49\textwidth}
        \centering
        \includegraphics[width=\linewidth]{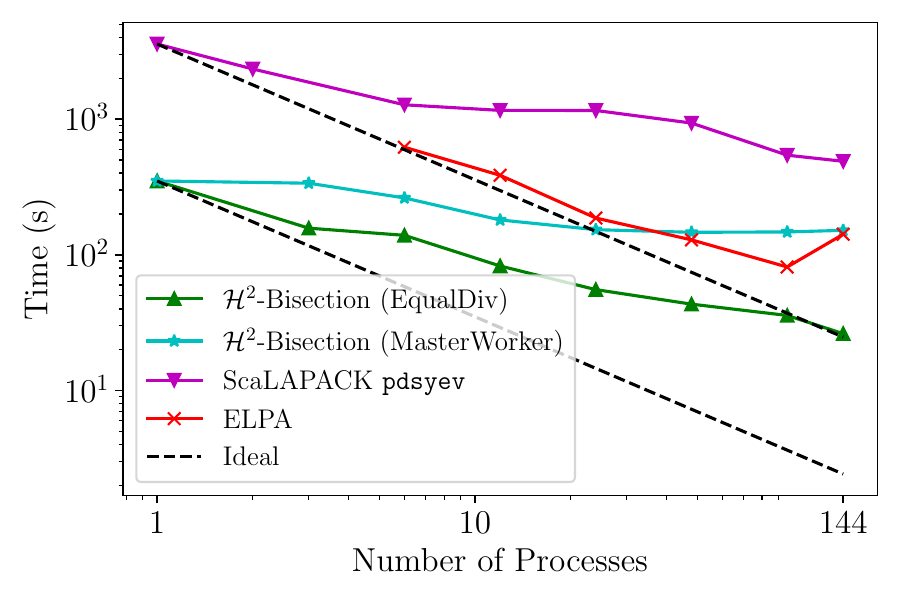}
        \caption{}
        \label{fig:laplace-par-eig-time}
    \end{subfigure}
    \hfill
    \begin{subfigure}[b]{0.49\textwidth}
        \centering
        \includegraphics[width=\linewidth]{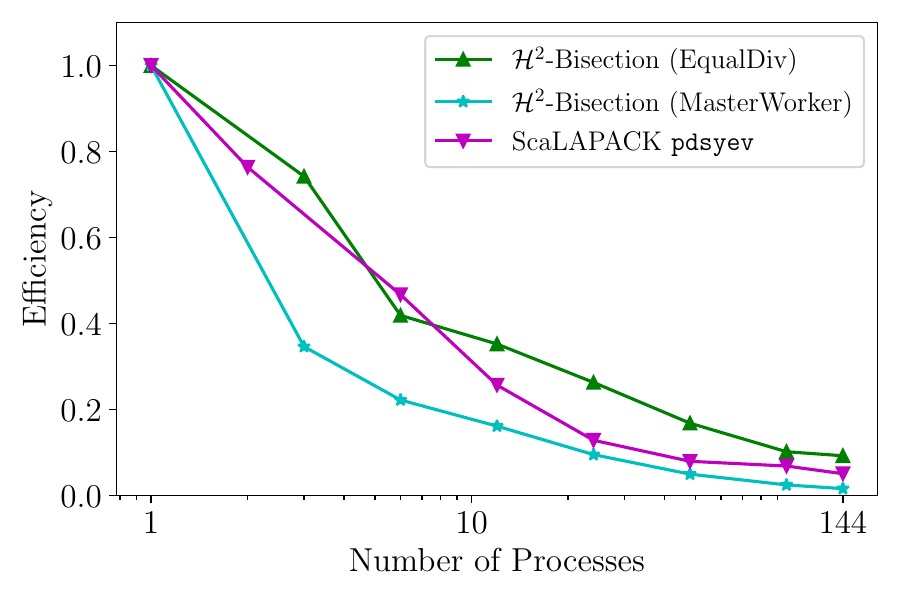}
        \caption{}
        \label{fig:laplace-par-eig-efficiency}
    \end{subfigure}
    \caption{Computing 100 eigenvalues of order 32,768 synthetic Laplace matrix with $\epsilon_{ev}=10^{-5}$ using up to 3 compute nodes. (a) Computation time. (b) Process efficiency.}
    \label{fig:laplace-par-eig-res}
\end{figure}
Finally, we show the scalability of the parallel algorithms discussed in Section~\ref{subsec:parallelization} in computing 100 eigenvalues of a given $n \times n$ matrix where $k_0 = \lfloor n/2 \rfloor - 50$ and $k_1 = k_0 + 99$. The starting interval of $[0,2048]$ was used for our bisection method. Small values of $m$ were chosen for the master-worker model to allow for load balancing. We used up to 144 MPI processes, where each process was mapped to one core. The subroutine {\ttfamily{MPI\_Dims\_create}} was used to determine the size of the 2D process grid for ScaLAPACK and ELPA.
Fig.~\ref{fig:laplace-par-eig-time} shows that our first parallel algorithm was faster than both ScaLAPACK and ELPA when using up to 144 processes. On the other hand, our second algorithm with the master-worker model did not scale as much, making it only outperform ScaLAPACK for this particular case. This is because the 100 target eigenvalues were closely clustered within a small interval so that the generation of tasks became a bottleneck and the workers spent most of their time waiting for task. The same reason is also related to the low efficiency shown in Fig.~\ref{fig:laplace-par-eig-efficiency}. When the target eigenvalues are clustered in an interval that is much smaller than the initial bisection interval, reusing the inertia from the computation of smallest and largest target eigenvalues will dramatically reduce the size of the starting interval for computing the subsequent eigenvalues, which naturally will also reduce their computation times. The smaller the interval containing the target eigenvalues than the initial bisection interval is, the more the speedup factor will decay. One way to address this is by supplying a sufficiently small initial bisection interval that contains all target eigenvalues.

\subsection{Computing the Target Eigenvalue for Electronic Structure Calculation}
\begin{figure}[hbtp]
    \centering
    \begin{subfigure}[b]{0.3\textwidth}
        \centering
        \includegraphics[width=0.8\linewidth]{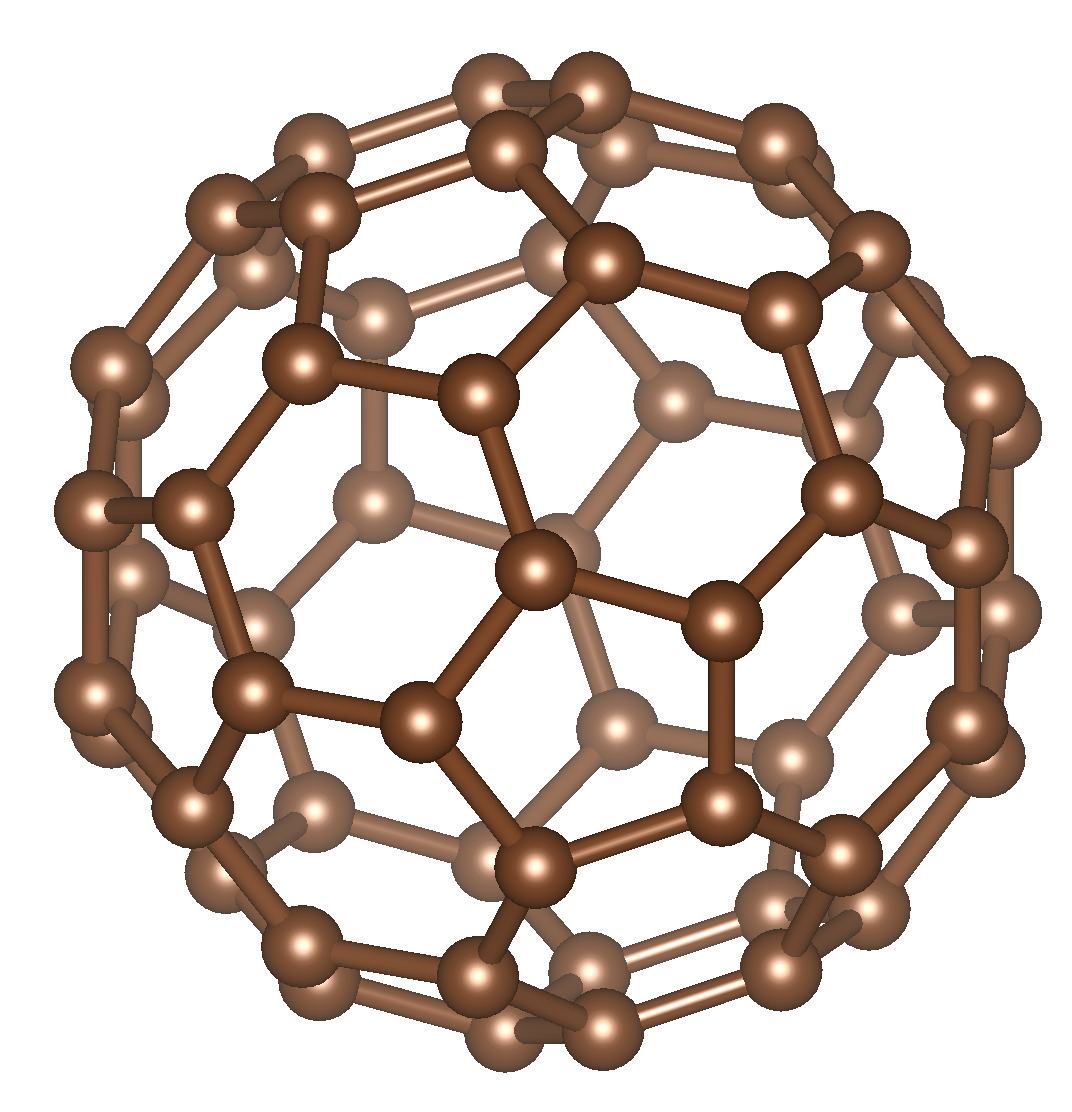}
        \caption{}
        \label{fig:fullerene-1}
    \end{subfigure}
    \hfill
    \begin{subfigure}[b]{0.69\textwidth}
        \centering
        \includegraphics[width=0.6\linewidth]{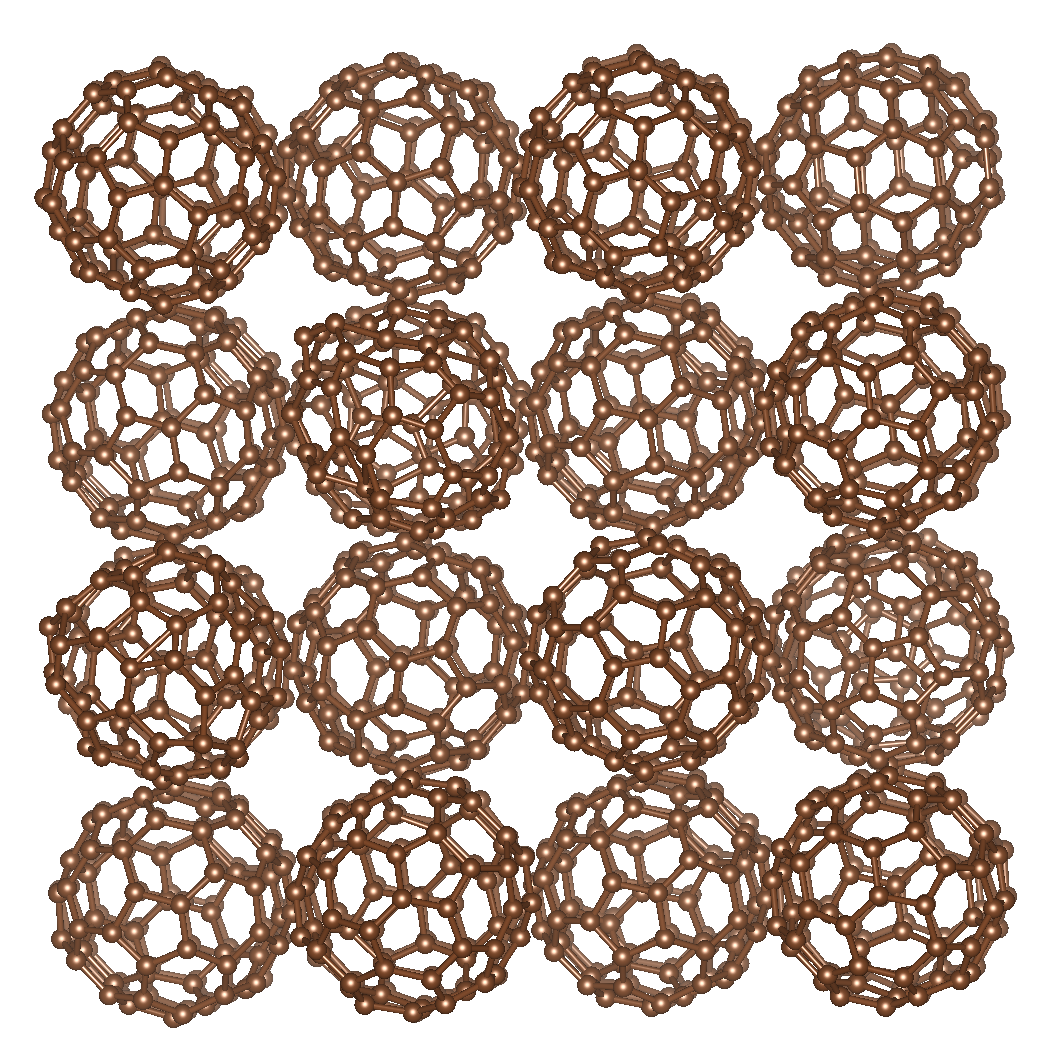}
        \caption{}
        \label{fig:fullerene-32}
    \end{subfigure}
    \caption{(a) A single fullerene mesh. (b) Arrangement of 32 fullerenes (2 identical layers of 4 $\times$ 4 fullerenes)}
    \label{fig:fullerene}
\end{figure}

Here we tested our method to solve the $k$-th eigenvalue problem arising from the electronic structure calculations of carbon nanomaterials composed of \textit{fullerene} (C\textsubscript{60}) allotropes~\cite{Omata2005_Fullerene}. A single fullerene takes the shape of a ball made up of 60 carbon atoms, as shown in Fig~\ref{fig:fullerene}. The real symmetric matrices corresponding to the standard eigenvalue problems were generated using ELSES quantum mechanical simulator~\cite{Hoshi2012_OrderNElectronicStructure}. The test materials that we used were formed by 64, 128, 256, and 512 fullerenes, which correspond to the matrix sizes of 15360, 30720, 61440, and 122880, respectively. We used the Hilbert space-filling curve to hierarchically partition the mesh where each fullerene ball formed the leaf block of the $\mathcal{H}^2$-matrix.

The target eigenvalue index $k$ is a material-specific value that is determined by the number of electrons in the material. In a typical case, the index is defined as $k=\lceil n/2 \rceil$ where $n$ is the matrix size. This target eigenvalue is important in determining the electronic properties of the material. For example, the difference between the $k$-th and ($k+1$)-th eigenvalue, which is referred to as the energy gap, can be used to determine the conductivity of metallic materials~\cite{Lee2018_KthEigenvalue}.

\begin{figure}[hbtp]
    \centering
    \begin{subfigure}[b]{0.49\textwidth}
        \centering
        \includegraphics[width=\linewidth]{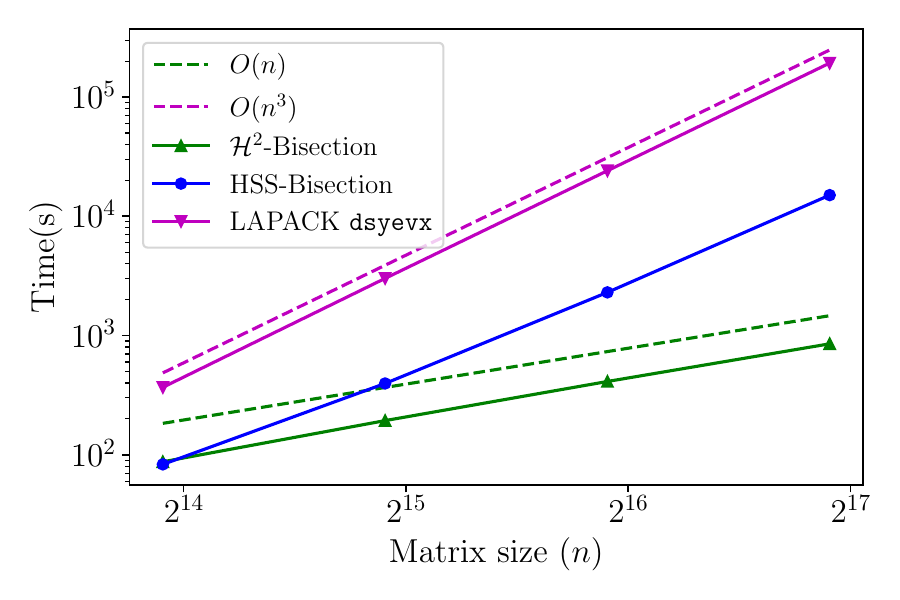}
        \caption{}
        \label{fig:elses-eig-time}
    \end{subfigure}
    \hfill
    \begin{subfigure}[b]{0.49\textwidth}
        \centering
        \includegraphics[width=\linewidth]{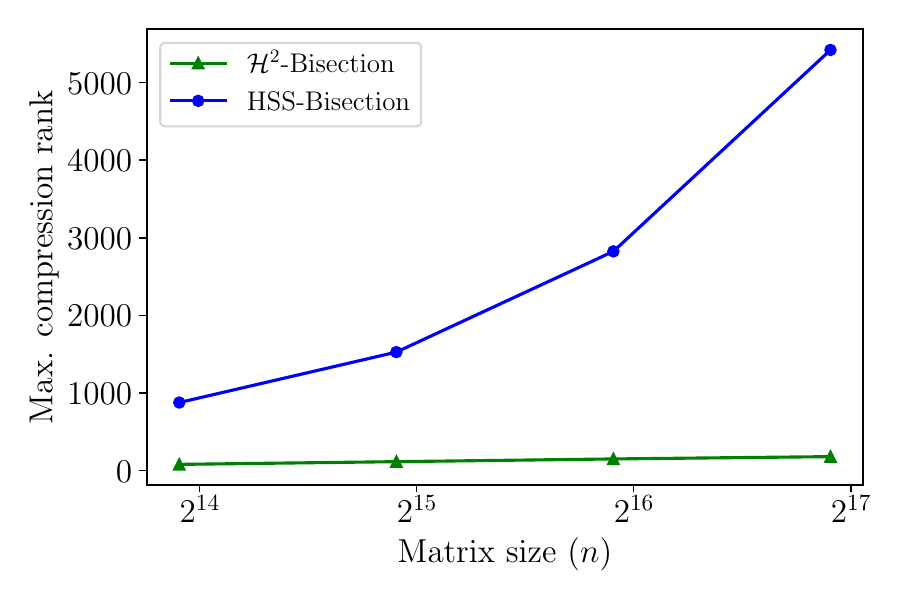}
        \caption{}
        \label{fig:elses-max-rank}
    \end{subfigure}
    \caption{Computing the target eigenvalue for electronic state calculations with $\epsilon_{ev}=10^{-3}$ using a single core of the compute node. (a) Computation time. (b) Maximum compression rank.}
    \label{fig:elses-result}
\end{figure}
Fig.~\ref{fig:elses-eig-time} shows the time to compute the target eigenvalue for each material with a known starting interval of $[-2,2]$ and accuracy of $\epsilon_{ev}=10^{-3}$ that is typically enough to produce an initial estimate for electronic structure calculations. Near linear scaling with respect to the matrix size can be observed from our $\mathcal{H}^2$ bisection method. On the other hand, the HSS bisection showed computation times that grow as $\mathcal{O}(n^3)$. This is because the matrices originate from 3D geometries where HSS could not maintain a constant rank as the problem size is increased, leading to a high compression rank that grew proportionally with the matrix size as shown in Fig~\ref{fig:elses-max-rank}. Whereas $\mathcal{H}^2$-matrix allows the off-diagonal blocks with high numerical rank to be subdivided further so that each compressed block could maintain a sufficiently small rank.

\section{Conclusion}
We have presented a linear time generalized LDL decomposition of $\mathcal{H}^2$ matrices and applied it to the bisection eigenvalue solver to compute the $k$-th eigenvalue efficiently with controllable accuracy. Numerical experiments showed that our method is orders of magnitude faster and require much less storage than the dense eigenvalue solver in LAPACK. Moreover, tests on electronic structure calculations of carbon nanomaterials demonstrated that our method outperforms the existing HSS-based bisection eigenvalue algorithm by achieving optimal complexity on matrices originating from 3D geometries.

Further, if more than one eigenvalue is desired, we can naturally parallelize the bisection method by exploiting the fact that slicing disjoint intervals are independent of each other and some of the computed inertia can be reused for the computation of other eigenvalues. This resulted in an efficient parallel algorithm that is faster than the existing state-of-the-art parallel dense eigenvalue solver such as ScaLAPACK and ELPA. It must be noted that although the computation time is scalable, the memory is not since each process has to hold an instance of the $\mathcal{H}^2$ matrix. Since the size of $\mathcal{H}^2$-matrix only grows linearly with the matrix dimension, this is still feasible for relatively large matrices. However, when the matrix does not fit on a single node anymore, lower-level parallelism has to be considered.

Future work concerns the parallelization of generalized LDL factorization to accelerate the computation of the $k$-th eigenvalue. The existing HSS-LDL factorization in~\cite{Xi2014_Bisection_HSS} can be extended to distributed memory systems and would also benefit from GPUs due to the inherent parallelism of the diagonal factorization at each level combined with the simple structure of HSS matrices. Furthermore, our $\mathcal{H}^2$-LDL factorization can also be parallelized by using a dynamic task-based execution model to fully utilize the parallelism among operations on the same as well as different levels of the factorization. In addition, the backward error bound of the $\mathcal{H}^2$-LDL factorization should be addressed too in order to improve its reliability when applied to the bisection eigenvalue algorithm.
\section*{Acknowledgements}
We sincerely thank Prof. Takeo Hoshi (Department of Applied Mathematics and Physics, Tottori University) for providing us with the models and matrices data corresponding to the electronic structure calculations of carbon nanomaterials.
This work was supported by JSPS KAKENHI Grant Number JP20K20624, JP21H03447, JP22H03598, JP23H00490. This work is supported by ”Joint Usage/Research Center for Interdisciplinary Large-scale Information Infrastructures” in Japan (Project ID: jh230009-NAHI).

\bibliographystyle{plain}
\bibliography{main}

\end{document}